\documentclass[final]{siamltex}
\usepackage{amsmath,amssymb}

\newtheorem{rem}[theorem]{Remark}
\newtheorem{rems}[theorem]{Remarks}

\newtheorem{thm}{\bf{Theorem}}[section]
\newtheorem{cor}[thm]{\bf{Corollary}}
\newtheorem{prop}[thm]{\bf{Proposition}}
\newtheorem{lem}[thm]{\bf{Lemma}}
\newtheorem{defn}[thm]{Definition}

\newcommand{\ga}{\alpha}

\newcommand{\gb}{\beta}

\newcommand{\gC}{\Gamma}

\newcommand{\gd}{\delta}

\newcommand{\bbr}{\mathbb{R}}

\newcommand{\bbc}{\mathbb{C}}

\newcommand{\calf}{\mathcal{F}}

\newcommand{\call}{\mathcal{L}}

\newcommand{\calr}{\mathcal{R}}

\newcommand{\opc}{\operatorname{C}}

\newcommand{\opbuc}{\operatorname{BUC}}

\newcommand{\opb}{\operatorname{B}}

\newcommand{\opl}{\operatorname{L}}

\newcommand{\opid}{\operatorname{id}}

\newcommand{\opv}[2]{#1: \text{dom}(#1 )\subset #2 \longrightarrow #2}

\newcommand{\pd}{\partial}

\newcommand{\dhr}{\mathrel{\lhook\joinrel\relbar\kern-.8ex\joinrel

            \lhook\joinrel\rightarrow}}

\newcommand{\dom}{\operatorname{dom}}

\author{Patrick Guidotti, Department of Mathematics, University of California,
        Irvine}
\title{Optimal Regularity for a Class of Singular Abstract Parabolic Equations}

\begin{document}
\bibliographystyle{plain}

\maketitle

\bigskip
\begin{abstract}
A general class of singular abstract Cauchy problems is considered which naturally 
arises in applications to certain Free Boundary Problems. Existence of an associated 
evolution operator characterizing its solutions is established and is subsequently 
used to derive optimal regularity results. The latter are well known to be important 
basic tools needed to deal with corresponding nonlinear Cauchy Problems such as those 
associated to Free Boundary Problems. 
\end{abstract}

\vfill \noindent {\small {\bf Keywords} Singular Parabolic Initial Boundary Value Problems, 
Optimal Regularity, Evolution Operator.}

\section{Introduction} 
In this paper a class of abstract Cauchy problems in a Banach space $E_0$ 
$$
 \dot u -A(t)u=f(t)\, ,\: t>0\, ,
$$
is considered where the family of operators $A$ is allowed to have singular 
behavior in the origin. The paper will focus on the ``parabolic case'', that is, 
it will be assumed that the operators $A(t)$ generate analytic semigoups for 
fixed $t>0$. It will, however, not be assumed that the latter be strongly continuous. 
Conditions on the family $A$ (see \eqref{hypo1}-\eqref{hypo3}) are proposed which lead 
to the construction of an associated evolution operator 
$U_A$ (cf. \eqref{evop}). By means of it, maximal regularity results are obtained for the 
singular abstract Cauchy problem in the context of spaces of singularly H\"older continuous 
functions. A certain class of Free Boundary Problems \cite{Gui99} with initial onset of a phase is 
the prime motivation for the study presented here. The paper is organized as follows. In the next 
section the problem is put into context and the basic tools needed in its analysis are presented. 
In section \ref{evopcons} the evolution operator $U_A$ is constructed and, in the following section, 
it is used to prove maximal regularity results for the singular abstract Cauchy problem. 
In final section \ref{example} an example is considered of a parabolic problem in a 
space-time wedge and the connection to Free Boundary Problems is made more explicit.

\section{Preliminaries and Setting}
Let $E_0$ be a Banach space. An unbounded operator
$$
  \opv{A}{E_0}
$$
is called {\em sectorial} if it satisfies
\begin{align}\label{sectorial}
  \text{(i) } &\overline{\operatorname{dom}(A)}=E_0\, ,\: N(A)=\{0\}\, ,\:\overline{R(A)}=E_0\, ,\\
  \text{(ii) } &(0,\infty)\subset\rho(A)\text{ and }\| t(t-A)^{-1}\|_{\mathcal{L}(E_0)}
  \leq M\, ,\: t>0\, ,\: M>0\, . 
\end{align}
If $A$ only satisfies (ii) is usually called {\em pseudo-sectorial}. For any given pseudo-sectorial 
operator $A$ on $E_0$ there exists $\theta>0$ such that
\begin{equation}\label{angle}
  \rho(A)\supset\Sigma _{\theta}:=\{\lambda\in\opc\setminus\{0\}\, |\, |\arg(\lambda)|<\theta\}
  \text{ and } \sup_{\lambda\in\Sigma _\theta}\| \lambda(\lambda-A)^{-1}\|_{\mathcal{L}(E_0)}\leq c
\end{equation}
thus clarifying the name. If $\theta>\pi/2$, it is well-know that an {\em analytic semigroup} 
$T_A$ can be associated to any given pseudo-sectorial operator $A$ through the formula
\begin{equation}\label{dunford-schwarz}
  T_A(t)=e^{tA}=\frac{1}{2\pi i}\int_{\Gamma}\!e^{\lambda t}(\lambda -A)^{-1}\,d\lambda\, ,\: t>0\, ,
\end{equation}
where the path $\gC$ is given by
$$
  \Gamma=\{\lambda\in\bbc\, |\, \arg(\lambda)=\eta\, ,\: |\lambda |\geq r\}\cup
  \{\lambda\in\bbc\, |\, |\arg(\lambda)|\leq\eta\, ,\: |\lambda |=r\}
$$
and is oriented counterclockwise. The parameters $r$ and $\eta$ are chosen such that $r>0$ and 
$\eta\in(\pi/2,\theta)$. If the operator $A$ is densely defined, then the semigroup $T_A$ is 
{\em strongly continuous}, that is, it satisfies
$$
  \lim _{t\to 0+}T_A(t)x=x\, ,\: x\in E_0\, .
$$
Otherwise it is strongly continuous on $\overline{\operatorname{dom}(A)}$.

The most important properties of analytic semigroups (and an equivalent characterization) are
\begin{align}\label{holo_char1}
  \text{(i) }&e^{tA}(E_0)\subset\operatorname{dom}(A)\, ,\: t>0\\\label{holo_char2}
  \text{(ii) }& \|tAe^{tA}\| _{\mathcal{L}(E_0)}\leq c\, ,\: t\in[0,T]\, ,\: T>0\, .
\end{align}
The semigroup $T_A$ is called exponentially decaying if it satisfies
\begin{equation}\label{exp_decaying}
  \| T_A(t)\| _{\mathcal{L}(E_0)}\leq ce^{-\omega t}\, ,\: t>0
\end{equation}
for some $\omega>0$. The collection of generators $A$ of analytic semigroups satisfying \eqref{exp_decaying} 
for some $c$ is denoted by $\mathcal{H}^-(E_0,\omega)$. Semigroups are useful in the analysis of 
{\em abstract Cauchy Problems} (ACP)
\begin{equation}\label{acp}
  \dot u-Au=f(t)\, ,\: u(0)=x\, .
\end{equation}
If $f\in\opl _1(0,T;E_0)$ and $x\in\ E_0$, a function $u\in\opc\bigl([0,T],E_0\bigr)$ satisfying 
\begin{equation}\label{vcf}
  u(t)=T_A(t)x+\int _0^t T_A(t-\tau)f(\tau)\, d\tau\, ,\: t\in[0,T]\, ,
\end{equation}
is called {\em mild solution} of \eqref{acp}. If the abstract Cauchy Problem is non-autonomous, 
that is, if $A$ depends on the time variable, then mild solutions of
\begin{equation}\label{nacp}
  \dot u-A(t)u=f(t)\, ,\: u(0)=x
\end{equation}
are given by
\begin{equation}\label{evop_vcf}
  u(t)=U_A(t,0)x+\int _0^t U_A(t,\tau)f(\tau)\, d\tau\, ,
\end{equation}
if it can be shown that an {\em evolution operator} $U_A$ associated to the family $A$ exists. 
The latter is a two-parameter family
\begin{equation}\label{evopdef}
  \{ U_A(t,\tau)\, |\, t\in[0,T]\, ,\:\tau\in[0,t]\}\subset\mathcal{L}(E_0)
\end{equation}
satisfying
\begin{align*}
  \text{(i) }& U_A(t,t)=\opid _{E_0}\, ,\: U_A(t,\tau)U_A(\tau,s)=U_A(t,s)\, ,\: 
  0\leq s\leq\tau\leq t\leq T\, ,\\
  \text{(ii) }& \partial _t U_A(t,\tau)x=A(t)U_A(t,\tau)x\, ,\: \tau <t\leq T\, ,\: x\in E_0\, .
\end{align*}
Classical results show that such an evolution operator exists on some regularity assumptions on 
the family $A$, usually of H\"older type. In \cite{Tan60,Sob66,Ama95}, in the case of densely 
defined family $A$, it is obtained as the solution to the weakly singular Volterra integral equation
\begin{equation}\label{vie}
  U_A(t,s)=e^{(t-s)A(s)}-\int _s^tU_A(t,\tau)[A(\tau)-A(s)]e^{(\tau-s)A(s)}\, d\tau\, .
\end{equation}
Another construction due to Da-Prato and Lunardi \cite{DaPG84,AT84,AT87,Lun95} is 
based on maximal regularity results for the autonomous abstract Cauchy problem combined with 
perturbation arguments but does not rely on the family $A$ having dense domains of definition.
They obtain 
$$
  U_A(t,s)x=W_A(t,s)x+e^{(t-s)A(s)}x
$$ 
from the solution $w=W_A(\cdot,s)x$ of
\begin{equation}\label{evopeq}
  \dot w(t)=A(s)w(t)+\bigl[ A(t)-A(s)\bigr]x\, ,\: w(s)=x\, .
\end{equation}
{\em Optimal or maximal regularity} results for \eqref{acp} or \eqref{nacp} can be described as 
follows: Find spaces $S$ and $E_0$ such that \eqref{acp} or \eqref{nacp} possess a unique solution 
$u\in S\bigl( (0,T],E_0\bigr)$ with
\begin{equation}\label{max-reg}
  \dot u\, ,\: Au\in S\bigl( (0,T],E_0\bigr)\text{ and } \|\dot u\| _{S\bigl( (0,T],E_0\bigr)}+
  \| Au\| _{S\bigl( (0,T],E_0\bigr)}\leq c\| f\| _{S\bigl( (0,T],E_0\bigr)}\, .
\end{equation}
The name obviously refers to the fact that $u$ enjoys as much regularity as is conceivably possible 
by virtue of it satisfying \eqref{acp} or \eqref{nacp}, respectively. It is known that restrictions 
apply to the choice of $S$ and $E_0$ for maximal regularity to hold. Counter-examples are given 
in \cite{Lun95} and \cite{LM99}.\\
Fractional powers and interpolation spaces play an important role in the theory of abstract Cauchy 
problems of parabolic type. 
Given a pseudo-sectorial operator $A$ it is always possible to define its {\em fractional powers} 
$(-A)^{\rho}$, $\rho >0$, as the inverses of the bounded operators
\begin{equation}\label{fractional_powers}
  (-A)^{-\rho}:=\frac{1}{\Gamma(\rho)}\int _0^\infty\, t^{\rho-1}e^{tA}\, dt\, .
\end{equation}
defined on their range, that is, with
$$
  \operatorname{dom}\bigl((-A)^{\rho}\bigr)=R\bigl((-A)^{-\rho}\bigr)
$$
For $\rho =1$, \eqref{fractional_powers} simply gives the resolvent of $A$ 
as the Laplace transform of the semigroup $e^{tA}$. The {\em interpolation spaces}
\begin{equation*}
  E^p_\alpha=D_A(\alpha ,p)\, ,\: \alpha\in(0,1)\, ,\: p\in[1,\infty]\text{ or }
  (\alpha ,p)=(1,\infty)\, ,
\end{equation*}
are defined by
\begin{align}\label{interpolation_space}
  E^p_\alpha :&=\big\{ x\in E_0\, \big |\, [t\mapsto v(t):=\| t^{1-\alpha -1/p}Ae^{tA}x\|]
  \in\opl _p(0,1)\big\}\\
  \| x\| _{E^p_\alpha}&=\| x\|+[x]_{\alpha ,p}:=\| x\|+\| v\|_{\opl _p(0,1)}\, ,
\end{align}
with the standard convention that $1/\infty =0$. The domain 
$\operatorname{dom}\bigl((-A)^{\alpha}\bigr)$ endowed with the graph norm $\| x\|_{D_\alpha}=
\| (-A)^{\alpha}x \|$  is denoted by $D_\alpha$ or $D((-A)^\alpha)$. It can be conveniently 
sandwiched between interpolation spaces
\begin{equation}\label{sandwich}
  E^1_\alpha\hookrightarrow D_\alpha\hookrightarrow E^\infty _\alpha\, ,\: \alpha\in(0,1)
\end{equation}
and satisfies the interpolation property
$$
  \|(-A)^{\alpha}x \|\leq c\,\| Ax\|^\alpha\| x\| ^{1-\alpha}\, .
$$
One of the main reasons to consider maximal regularity results is dealing with 
fully nonlinear Cauchy problems like
\begin{equation}\label{fnacp}
  \dot u=f(t,u)\, ,\: t>0\, ,\: u(0)=u_0\in E_0\, ,
\end{equation}
in $E_0$ or quasilinear problems like
\begin{equation}\label{qlacp}
  \dot u-A(t,u)u=f(t,u)\, ,\: t>0\, ,\: u(0)=u_0\in E_0\, ,
\end{equation}
The basic idea is to linearize \eqref{fnacp} (or \eqref{qlacp}) in $u_0$ and localize 
about $t=0$ to obtain
\begin{equation}\label{linearization}
  \dot u-A(0)u=G(u)\, ,\: t>0\, ,
\end{equation}
where $A(0)=DF(0,u_0)$ and $DG(0,u_0)=0$ and solve it by a fixed-point argument in a small 
time interval. Repeating the procedure the solution can be extended to its maximal interval 
of existence. This can only be done if maximal regularity results hold which are essentially 
equivalent to $ \partial_t -A(0)$ being invertible. In spite of the simplicity of this idea, its 
implementation is rather involved and needs the introduction of spaces of singular H\"older 
continuous functions in order to successfully deal with the kernel singularity \eqref{holo_char2} 
in the origin. A comprehensive exposition of this theory is given in \cite{Lun95}, where the 
author also presents a variety of examples which clearly attest to its wide range of applicability. 
One may also consult \cite{A90}.

The focus of this paper is on {\em singular} abstract Cauchy problems like
\begin{equation}\label{sacp}
  \dot u-A(t)u=f(t)\, ,\: t>0\, ,
\end{equation}
in a Banach space $E_0$. The family $A$ is allowed to behave singularly in the origin. In this 
case \eqref{sacp} is the appropriate model problem since there is no autonomous counterpart 
to speak of. One of the main goals of the paper is to find abstract but flexible conditions on 
the singular family $A$ which allow for the construction of an associated evolution operator and 
eventually lead to maximal regularity results for \eqref{sacp}. Again the prime applications would 
be fully nonlinear counterparts of \eqref{sacp}, which can be written as
\begin{equation}\label{nsacp}
  \dot u-A(t,u)u=f(t,u)\, ,\: t>0\, .
\end{equation}
In this case the problem needs not only to be linearized as in the regular case but also to be 
expanded in the singularity in order to capture its leading order behavior and perturbe around it.
A notation closer to \eqref{qlacp} rather than \eqref{fnacp} is used to draw reader's attention 
to the singular behavior in the origin of the generator families $A(\cdot,u)$ involved. 
A class of Free Boundary Problems, considered in \cite{Gui99} for the first time in higher dimensions 
and by many authors in a one dimensional setting \cite{AGL66,AS78,FP77a,FP77b,FP77c,FMP86,COEr88,
Gui96c}, leads naturally to equations of singular type and motivates the study presented here.

Singular families of type
\begin{equation}\label{oldfam}
 A(t)=\frac{A}{t^k}\, ,\: k>1\, ,
\end{equation}
and natural generalizations thereof were considered in \cite{Gui96a} where an associated 
evolution operator was constructed by using \eqref{vie} and results for a corresponding 
class of quasilinear equations were obtained. In the purely linear case, they were also considered 
by \cite{Fa85,W98} were operator sum/product techniques were used in a $\opl _p$-setting. 

Here the aim is to extend those results to obtain maximal regularity in spaces of (singular) 
H\"older continuous functions for a non-densely defined family of generators and for a 
generalization of \eqref{oldfam} in which the singularity is allowed to affect the operator 
in an anisotropic way. A simple example is given by
\begin{equation}\label{proto_ex}
  A(t)=B+\frac{C}{t^k}
\end{equation}
where $B,C$ are given generators of analytic semigroups in the sense explained of 
\eqref{dunford-schwarz}. In applications the operators $B$ and $C$ are usually differential 
operators acting on distinct spacial variables. However, they typically do not commute as 
they in general have non-constant coefficients. All the results mentioned above do not apply 
to the setting of this paper and only those in \cite{Gui96a} could be adapted to it but 
merely in the event that $B$ and $C$ were commuting operators. This latter case is not of 
much interest since it never occurs in practical applications.

Moreover the conditions previously used in \cite{Gui96a} are ad-hoc in that they rely on 
the singularity being of a given explicit type (power type) in order to construct the associated 
evolution operator. In \cite{Fa85,W98} no evolution operator is constructed but the singularity 
is assumed to be given by a simple function of time.

Here a more abstract condition is obtained which seems quite natural in the construction of the 
evolution operator as the proofs will show and which allows for a wider class of singularities 
including all those previously considered in the mentioned papers. These conditions read
\begin{align}\label{hypo1}
  \text{(i) }&A(t)\in\mathcal{H}^-(E_0,\omega)\, ,\: t>0\, ,\\\notag
  \text{(ii) }& \|\bigl[ A(t)-A(s)\bigr]A^{-1}(\tau)\|_{\mathcal{L}(E_0)}\leq c\frac{t-s}{t}
  \text{ and }\\\label{hypo2} &\|\bigl[ A(t)-A(s)\bigr](-A)^{-\rho}(\tau)\|_{\mathcal{L}(E_0)}
  \leq c(t-s)\, ,\\\label{hypo3}
  \text{(iii) }&\lim _{t\to 0}A^{-1}(t)=0\, ,
\end{align}
for some $\rho\in(1,2)$ and $0<\tau\leq s\leq t\leq T$. As an example one can consider 
$A(t)=B+\frac{C}{t^2}$ satisfying \eqref{hypo1} and such that $C$ is an invertible pseudo-sectorial 
operator, then \eqref{hypo2}-\eqref{hypo3} are easily seen to be satisfied.

For the sake of completeness we observe that maximal regularity results have been obtained, 
in the regular case not considered here, in a variety of other function spaces other than 
singular H\"older spaces and with the base space $E_0$ substituted in particular by interpolation 
spaces of type \eqref{interpolation_space}. The interested reader may consult 
\cite{DaPG84,Lun95,Ama95,DV87,KW04} for instance. Some of these results cannot be obtained for 
the singular (anisotropic case), whereas others do not fit the requirements imposed by the 
applications to Free Boundary Problems. In spite of the fact that they would be of theoretical and 
possibly of practical interest, they are not considered in this paper.
\section{Construction of the evolution operator}\label{evopcons}
In order to construct the evolution operator $U$ associated to a singular family satisfying 
\eqref{hypo1}-\eqref{hypo3} it is necessary to work in spaces of singular H\"older 
continuous functions.
\begin{defn}
Let $\alpha\in(0,1)$ and $T>0$ be given and let $E$ be a Banach space.
\begin{equation}\label{sing_hoelder_space}
 \opc ^\ga_\ga\bigl( (0,T],E\bigr):=\Big\{ v\in\opb\bigl( (0,T],E\bigr)\, \Big |\, 
 [t\mapsto t^\ga v(t)]\in\opc ^\ga\bigl( (0,T],E\bigr)\Big\}\, .
\end{equation}
Endowed with the norm given by
$$
 \| v\| _{\alpha,s}:=\| v\| _\infty + [v]_{\ga,s}
$$
this space becomes a Banach space. Hereby we denoted by $[\cdot ]_{\ga,s}$ the weighted H\"older 
seminorm
$$
 [v]_{\ga,s}:=[(\cdot)^\ga v]_\ga :=\sup_{0<t\neq s\leq T}\frac{\| t^\alpha v(t)-s^\alpha v(s)\| 
 _E}{|t-s|^\ga}
$$
where $[\cdot]_\ga$ denotes the regular H\"older seminorm. We shall also make use of the space
\begin{equation}\label{0_hoelder_space}
 \opc ^\ga_0\bigl( (0,T],E\bigr):=\Big\{ v\in\opc ^\ga( [0,T],E\bigr)\, \Big |\, v(0)=0\Big\}\, .
\end{equation}
For $\beta\in(0,1)$ the additional space
\begin{multline}\label{gen_sing_hoelder_space}
 \opc ^\ga_{\alpha,\beta}\bigl( (0,T],E\bigr):=\Big\{ v:(0,T]\to E\, \Big |\\ [t\mapsto t^\gb v(t)]
 \in\opb\bigl((0,T],E\bigr)\, ,\: [t\mapsto t^{\alpha+\beta} v(t)]\in\opc ^\ga\bigl( (0,T],E\bigr)
 \Big\}\, .
\end{multline}
is also defined and endowed with its natural norm 
$$
  \|v\| _{\alpha,\beta}=\|(\cdot)^\beta v\|_\infty +[(\cdot)^{\alpha+\beta} v]_\alpha\, .
$$
The space
$$
  \opb _{\beta}\bigl((0,T],E\bigr):=\big\{ v:(0,T]\to E\,\big |\, [t\mapsto t^\gb v(t)]
  \in\opb\bigl((0,T],E\bigr)\big\}
$$
will also be useful.
\end{defn}\\
Fix $x\in E_0$. Then it is natural to look for $U(\cdot ,s)x$ as the solution of
\begin{equation}\label{evopeq1}
  \dot u=A(t)u\, ,\: t\in(s,T]\, ,\: u(s)=x\, ,
\end{equation}
which is rewritten as
\begin{equation}\label{evopeq2}
  \dot w=A(t)w+\bigl[ A(t)-A(s)\bigr]e^{(t-s)A(s)}x\, ,\: t\in(s,T]\, ,\: u(s)=0\, .
\end{equation}
by setting $w(t):=u(t)-e^{(t-s)A(s)}x$. If \eqref{evopeq2} can be solved and denoting by 
$W(\cdot,s)x$ its solution in that case, the evolution operator is then simply given by
\begin{equation}\label{evop}
  U(t,s)=W(t,s)+e^{(t-s)A(s)}\, .
\end{equation}
The next theorem establishes existence for \eqref{evopeq2} in the space 
\eqref{gen_sing_hoelder_space}.
\begin{thm}\label{evopex}
Assume that $A$ satisfies \eqref{hypo1}-\eqref{hypo3} with $\rho\in(1,2)$ and let 
$f\in\opc ^\ga_{\alpha,\rho -1}\bigl( (s,T],E\bigr)$ for $\alpha\in(0,1)$, $s\in(0,T)$. 
Then the solution $w$ of
\begin{equation*}
  \dot w=A(t)w+\underset{=:g_s(t)}{\underbrace{\bigl[ A(t)-A(s)\bigr]e^{(t-s)A(s)}x}}
  +f(t)\, ,\: t\in (s,T]\, ,\: w(0)=0\, ,
\end{equation*}
satisfies
\begin{gather}\label{max_reg_0_est}
  w,\dot w,A(s)w\in\opc ^\ga_{\alpha,\rho -1}\bigl( (s,T],E_0\bigr)\, ,\: \dot w\in
  \opb _{\alpha+\rho -1}\bigl( (s,T],E^\infty_\alpha\bigr)\, ,\\
  \| Aw\| _{\opc ^\ga_{\alpha,\rho -1}E_0} +\| \dot w\| _{\opc ^\ga_{\alpha,\rho -1}E_0} +
  \| \dot w\| _{\opb _{\alpha+\rho -1}E^\infty _\alpha}\leq c\bigl( \| x\| +
  \| f\| _{\opc ^\ga_{\alpha,\rho -1}E_0}\bigr)\, .
\end{gather}
\end{thm}
\begin{proof}
The solution is constructed as the unique fixed-point of $\Phi$ in the space 
$\opc ^\ga_{\alpha,\rho -1}\bigl( (s,T],D(A(s))\bigr)$ where $\Phi(v)$ is defined as the solution of
\begin{equation*}
  \dot w=A(s)w+\bigl[ A(t)-A(s)\bigr]\bigl[ v+e^{(t-s)A(s)}x \bigr]+f(t)\, ,\: t\in (s,T]
  \, ,\: w(0)=0\, .
\end{equation*}
\underline{Step 1}: First it is checked that $g_s\in\opc ^\ga_{\alpha,\rho -1}\bigl( (s,T],
E_0\bigr)$. In fact
\begin{multline*}
  \| \bigl[ A(t)-A(s)\bigr]e^{(t-s)A(s)}x\|=\\\|\bigl[ A(t)-A(s)\bigr]\bigl(-A(s)\bigr)^{-\rho}
  \bigl( -A(s)\bigr)^\rho e^{(t-s)A(s)}x\|\leq c\frac{1}{(t-s)^{\rho -1}}\, ,
\end{multline*}
which gives $g_s\in\opb _{\rho -1}\bigl( (s,T],E_0\bigr)$. Next take $0<s<s+\varepsilon
\leq r\leq t\leq T$ and consider
\begin{multline*}
  \varepsilon^{\alpha+\rho -1}\|g_s(t)-g_s(r)\|\leq\varepsilon^{\alpha+\rho -1}
  \|\bigl[ A(t)-A(r)\bigr]\bigl(-A(s)\bigr)^{-\rho}\bigl( -A(s)\bigr)^\rho e^{(t-s)A(s)}x\| 
  \\+\varepsilon^{\alpha+\rho -1}\|\|\bigl[ A(r)-A(s)\bigr]\int_{r-s}^{t-s}\!A(s)e^{\sigma A(s)}
  x\,d\sigma\|\\\leq c\varepsilon^{\alpha+\rho -1}(t-r)\frac{1}{(t-s)^\rho}\| x\|+c\varepsilon
  ^{\alpha+\rho -1}(r-s)\int _{r-s}^{t-s}\frac{1}{\sigma ^{1+\rho}}\, d\sigma\\\leq
  c\| x\|\bigl[ (t-r)^\alpha(t-r)^{1-\alpha}\frac{\varepsilon^{\alpha+\rho -1}}{t-s}+
  \varepsilon^{\alpha+\rho -1}(r-s)\frac{(t-s)^\rho-(r-s)^\rho}{(r-s)^\rho(t-s)^\rho}\bigr]
  \\\leq c\| x\| (t-r)^\alpha\, ,
\end{multline*}
since $(t-s)^\rho-(r-s)^\rho\leq c(t-s)^{\rho -1}(t-r)$. It is therefore seen that
$$
  g_s\in\opc ^\alpha_{\alpha+\rho -1}\bigl( (s,T],E_0\bigr)\text{ and }\| g_s\| _{\alpha,\rho -1}
  \leq c\| x\|\, ,
$$
for a constant $c$ which does not depend on $s$.\\
\underline{Step 2}: Next it is shown that $\Phi$ is a contractive self-map on 
$\opc ^\ga_{\alpha,\rho -1}\bigl( (s,s+\delta],D(A(s))\bigr)$ provided $\delta>0$ is 
small enough. In order to do so, it is enough to show contractivity together with 
\begin{equation}\label{phi-reg}
  \bigl[ A(\cdot)-A(s)\bigr]v\in\opc ^\alpha _{\alpha+\rho -1}\bigl( (s,T],E_0\bigr)
\end{equation}
because the existence and regularity of the solution $\Phi(v)$ then follows from 
known maximal regularity results for the regular case (cf. \cite[Theorem 4.3.7]{Lun95}).
In order to show \eqref{phi-reg} observe first that
\begin{multline*}
  \| (t-s)^{\rho -1} \bigl[ A(t)-A(s)\bigr]\bigl(-A(s)\bigr)^{-1}A(s)v(t)\|\\
  \leq\frac{t-s}{t}\| (t-s)^{\rho -1}A(s)v(t)\|\leq\frac{\delta}{s}\| v\| _{\opb _{\rho -1}
  D( A(s))}\, .
\end{multline*}
Next take $0<s<s+\varepsilon\leq r\leq t\leq T$ and consider
\begin{multline*}
  \varepsilon^{\alpha+\rho -1}\|\bigl[ A(t)-A(s)\bigr]v(t)-\bigl[ A(r)-A(s)\bigr]v(r)\|\\
  \leq\varepsilon^{\alpha+\rho -1}\|\bigl[ A(t)-A(r)\bigr]v(t)\|+\varepsilon^{\alpha+\rho -1}
  \|\bigl[ A(r)-A(s)\bigr](v(t)-v(r))\|\\\leq \varepsilon^{\alpha+\rho -1}\frac{t-r}{t}
  \frac{1}{\varepsilon^{\rho -1}}\|v\| _{\opb _{\rho -1}D( A(s))}+\frac{r-s}{r}(t-r)^\alpha
  \| v\| _{\opc ^\alpha _{\alpha+\rho -1}D( A(s))}\\\leq c\frac{\delta}{s}\| v\| _{\opc ^\alpha 
  _{\alpha+\rho -1}D( A(s))}\, ,
\end{multline*}
for a constant $c$ independent of $s$. This gives \eqref{phi-reg} and shows that
\begin{equation*}
  \|\bigl[ A(\cdot)-A(s)\bigr]v\| _{\alpha ,\rho -1,E_0}\leq c\frac{\delta}{s}\| 
  v\|_{\alpha ,\rho -1, D(A(s))}\, .
\end{equation*}
on the interval $(s,s+\delta]$. Latter estimate also gives that
\begin{equation*}
 \|\Phi(v_1)-\Phi(v_2)\| _{\alpha,\rho -1,D(A(s))}\leq c\frac{\delta}{s}
 \| v_1-v_2\| _{\alpha,\rho -1,D(A(s))}
\end{equation*}
in view of \cite[Theorem 4.3.7]{Lun95} and the linearity of the equation. It is then 
clear that $\Phi$ is a contraction for $\delta<<1$ and it is easy to obtain the inequality
$$
  \| w\| _{\alpha,\rho -1,D(A(s))}\leq c\| g_s\| _{\alpha,\rho -1,E_0}\leq c\| x\|\, .
$$
for the unique solution $w$ for a constant $c$ independent of $s$. By using further 
results from the regular theory \cite[Proposition 6.1.3]{Lun95} the solution can be 
continued to the full interval maintaining the inequality.
\end{proof}

It follows that $W$ (for $f\equiv 0$) satisfies the slightly better inequality
\begin{equation}\label{evopestpar}
  \| A(\tau)W(t,\tau)\| _{\mathcal{L}(E_0)}\leq c\frac{1}{(t-\tau)^{\rho-1}}\, ,\: 
  t\in(s,T]\, .
\end{equation}
Taking \eqref{evop} into account it is concluded that
\begin{equation}\label{evopestfull}
  \| A(\tau)U(t,\tau)\| _{\mathcal{L}(E_0)}\leq c\frac{1}{(t-\tau)}\, ,\: 
  t\in(s,T]\, .
\end{equation}
It is important to point out that the constant appearing in both estimates is independent 
of $\tau$ as follows from the proof of theorem \ref{evopex}.
\begin{rem}
In the regular case there is no restriction in the choice of the exponent $\rho -1$, whereas 
here it is determined by the singularity through \eqref{hypo2}.
\end{rem}
\begin{cor}\label{evopin0}
Let $A$ satisfy \eqref{hypo1}-\eqref{hypo3}. Then there is a unique evolution operator $U$ 
associated to $A$ defined for $T\geq t\geq\tau >0$. It can be extended to $\tau =0$ by 
setting
$$
  U(t,0)=0\, ,\: 0<t\leq T\, .
$$ 
\end{cor}
\begin{proof}
The claim follows from
\begin{multline*}
  \| U(t,\tau)\| _{\mathcal{L}(E_0)}\leq \| A^{-1}(\tau)\| _{\mathcal{L}(E_0)}
  \| A(\tau)U(t,\tau)\| _{\mathcal{L}(E_0)}\\\leq c\frac{1}{t-\tau}\| A^{-1}(\tau)
  \| _{\mathcal{L}(E_0)}\longrightarrow 0\: (\tau\to 0)  
\end{multline*}
in view of \eqref{hypo3}.
\end{proof}\\
The evolution operator allows to characterize solutions of \eqref{sacp} through the 
classical variation-of-constant-formula.
\begin{prop}\label{svcf}
Lt $f\in\opl _1\bigl((0,T),E_0\bigr)$. Then any bounded mild solution $u$ of \eqref{sacp} 
is given by
$$
  u(t)=\int _0^t U(t,\tau)f(\tau)\, d\tau\, .
$$
\end{prop}
\begin{proof}
The assumption on the function $f$ together with ensures the existence of the variation 
of constant integral. In the event that the solution $u$ to \eqref{sacp} exists and is 
bounded, it coincides with the solution 
$u_\gd$ of \eqref{sacp} on $t\geq\gd$ corresponding to the initial condition $u_\delta(\delta)=
u(\delta)$. It is then necessarily given by
$$
 u(t)=u_\gd(t)=U(t,\delta)u(\gd)+\int _\gd ^t U(t,\tau)f(\tau)\, d\tau\, ,\: t\geq\gd\, .
$$
Since the convolution integral exists for $\gd =0$ the second term converges to \eqref{svcf}. By 
assumption the first can be estimated as follows
$$
 \| U(t,\gd)u(\gd)\| _{ E_0}\longrightarrow 0\, (\gd\to 0\, ,\: t>0)
$$
by virtue of the solution's boundedness and corollary \ref{evopin0}.
The function $u$ is therefore a mild solution of \eqref{sacp} for $t>0$.
\end{proof}
\begin{rems}
{\bf (a)} It now becomes clear why \eqref{sacp} is formulated without any initial condition. 
Bounded solutions naturally emanate from 0.\\
{\bf (b)} If additional information is available about the rate of vanishing of 
$A^{-1}(t)$ in $t=0$, it is possible to weaken the assumption on $f$ to 
$$
  \bigl[ t\to t^pf(t)\bigr]\in\opl _1\bigl((0,T),E_0\bigr)
$$ 
for an appropriate power $p>1$.
\end{rems}
\section{Maximal regularity}\label{maxreg}
In order to prove maximal regularity results for \eqref{sacp} a couple of lemmata are needed.
\begin{lem}\label{lunardi}
Let $x\in E_0$ and $A$ be the generator of a not necessarily strongly continuous analytic 
semigroup. Then $\int _0^te^{sA}x\, ds\in D(A)$ and
$$
  A \int _0^te^{sA}x=e^{tA}x-x\, .
$$
\end{lem}
\begin{proof}
The proof would be completely obvious if the semigroup were strongly continuous. On the 
given assumptions it needs a little more care but a proof can be found in \cite{Lun95}.
\end{proof}
\begin{lem}\label{intest}
Assume that $A$ satisfies assumptions \eqref{hypo1}-\eqref{hypo3}. Then
\begin{equation*}
  \| A(t)\bigl[e^{(t-\tau)A(\tau)}-e^{(t-\tau)A(t)}\bigr]\| _{\mathcal{L}(E_0)}
  \leq c\frac{1}{t}\, ,\: 0<\tau< t\leq T\, .
\end{equation*}
\end{lem}
\begin{proof}
Since
$$
  e^{(t-\tau)A(\tau)}=\frac{1}{2\pi i}\int _\Gamma e^{\lambda(t-\tau)}(\lambda -A(\tau))^{-1}
  \, d\lambda
$$
the estimand can be rewritten as
\begin{multline*}
  -\frac{1}{2\pi i}\int _\Gamma e^{\lambda(t-\tau)} A(t)(\lambda -A(t))^{-1}\bigl[
  A(\tau)-A(t)\bigr]A^{-1}(\tau)A(\tau)(\lambda -A(\tau))^{-1}\, d\lambda\, .
\end{multline*}
Assumptions \eqref{hypo1} and \eqref{hypo2} give
$$
 \| A(t)(\lambda -A(t))^{-1}\bigl[ A(\tau)-A(t)\bigr]A^{-1}(\tau)A(\tau)
 (\lambda -A(\tau))^{-1}\| _{\mathcal{L}(E_0)}\leq c\frac{t-\tau}{t}
$$
and the claim follows by direct estimation of the integral.
\end{proof}
\begin{thm}\label{max_reg_0}
Assume that $A$ satisfies assumptions \eqref{hypo1}-\eqref{hypo3} and let 
$f\in\opc ^\alpha _0\bigl((0,T],E_0\bigr)$ for some $\alpha\in(0,1)$. 
Then the solution $u$ of \eqref{sacp} given by \eqref{svcf} on $(0,T]$ satisfies
$$
  \dot u ,Au\in\opc ^\alpha _0\text{ and } \| \dot u\| _{\ga}+\| Au\| _{\alpha}\leq 
  c\| f\| _\alpha\, .
$$
\end{thm}
\begin{proof}
\underline{Step 1}: First H\"older continuity in the origin is shown. To that end it is 
enough to show that
$$
  \| A(t)\int _0^t U(t,\tau)f(\tau)\, d\tau\|\leq c t^\alpha
$$
since then the equation gives the corresponding estimate for $\dot u(t)$. Using the 
decomposition \eqref{evop} of the evolution operator it is seen that
\begin{multline*}
  \int _0^t U(t,\tau)f(\tau)\, d\tau=\int _0^t W(t,\tau)f(\tau)\, d\tau +\int _0^t 
  e^{(t-\tau)A(\tau)}[f(\tau)-f(t)]\, d\tau\\+\int _0^t\bigl[ e^{(t-\tau)A(\tau)}-e^{(t-\tau)A(t)}
  \bigr]f(t)\, d\tau+\int _0^t e^{(t-\tau)A(t)}f(t)\, d\tau\\ I+II+II+IV\, .
\end{multline*}
The various terms can then be estimated as follows
\begin{multline*}
  \| A(t)I\|\leq c\int _0^t\frac{\tau ^\alpha}{t-\tau}\, d\tau\leq ct^\alpha
  \text{ by \eqref{evopestpar}}\\
  \| A(t)II\|\leq c\int _0^t\frac{1}{(t-\tau)^{1-\alpha}}\, d\tau=ct^\alpha
  \text{ by assumption}\\
  \| A(t)III\|\leq c\int _0^t \frac{1}{t^{1-\alpha}}\, d\tau=ct^\alpha
  \text{ by lemma \ref{intest}}
\end{multline*}
As for the last one has $A(t)IV=e^{tA(t)}f(t)-f(t)$ by lemma \ref{lunardi} which 
gives the desired estimate by the assumptions on $f$.\\
\underline{Step 2}: Away from the origin, solution properties should not deviate from 
the regular case. This is in fact confirmed by the following argument.
It follows from the previous step that
$$
  u(t)=\int_{0}^{t}\!U(t,\tau)\,d\tau\in D(A(t))\, ,\: t\in(0,T]\, .
$$
Thus $u(\delta)\in D(A(\delta))$ and, using theorem \ref{evopex} one gets that
$$
  \dot u(\delta)=A(\delta)u(\delta)+f(\delta)\in E^\alpha _\infty
$$
and that
\begin{multline*}
  \| u(\delta)\|_{D(A(\delta))}\leq c\| f\| _{\alpha,\rho -1}\leq c\| f\| _\alpha\\
  \text{ and }\|A(\delta)u(\delta)+f(\delta)\| _{E^\alpha _\infty}\leq c
  \| f\| _{\alpha,\rho -1}\leq c\| f\| _\alpha\, .
\end{multline*}
The embedding inequalities are a consequence of
\begin{multline*}
  \| t^{\rho -1}f(t)\|\leq ct^{\rho -1+\alpha}[f]_\alpha\text{ and of }\\
  \|\frac{t^{\alpha+\rho-1}f(t)-s^{\alpha+\rho-1}f(s)}{(t-s)^\alpha}\|\leq 
  \|\frac{t^{\alpha+\rho-1}-s^{\alpha+\rho-1}}{(t-s)^\alpha}f(t)\|+\|s^{\alpha+\rho-1}
  \frac{f(t)-f(s)}{(t-s)^\ga}\|\\\leq (t-s)^{\min(\alpha+\rho-1,1)}t^\alpha[f]_\alpha +
  s^{\alpha+\rho-1}[f]_\alpha\, .
\end{multline*}
Thus, using \cite[Theorem 4.3.1(iii)]{Lun95}, it follows that
$\dot u,Au\in\opc ^\alpha\bigl([\delta,T],E_0\bigr)$ for any given $\delta>0$ and
$$
  \|\dot u\| _\alpha+\|Au\| _\alpha\leq c\| f\| _\alpha
$$ 
which gives H\"older continuity everywhere away from the origin.
\end{proof}
\begin{thm}\label{max_reg_sing}
Assume that $A$ satisfies assumptions \eqref{hypo1}-\eqref{hypo3} and let 
$f\in\opc ^\alpha _\alpha\bigl((0,T],E_0\bigr)$ for some $\alpha\in(0,1)$. 
Then the solution $u$ of \eqref{sacp} given by \eqref{svcf} on $(0,T]$ satisfies
$$
  \dot u ,Au\in\opc ^\alpha _\alpha\text{ and } \| \dot u\| _{\alpha,s}+\| Au\| _{\alpha,s}
  \leq c\| f\| _{\alpha ,s}\, .
$$
\end{thm}
\begin{proof}
\underline{Step 1}: First consider regularity in the origin. From \eqref{evop} it follows that
\begin{multline*}
  A(t)\int_{0}^{t}\! U(t,\tau)f(\tau)\,d\tau=A(t)\int_{0}^{t}\! W(t,\tau)f(\tau)\,d\tau+
  \int_{0}^{t}\! e^{(t-\tau)A(\tau)}f(\tau)\,d\tau=I_1+I_2\, .
\end{multline*}
The first term leads to
$$
 \| I_1\|\leq \int_{0}^{t}\!\frac{1}{(t-\tau)^{\rho-1}}\,d\tau=ct^{2-\rho}\, .
$$
$I_2$ needs to be further split
\begin{multline*}
  I_2=A(t)\int_{0}^{t}\! U(t,\tau)[f(\tau)-f(t)]\,d\tau+A(t)\int_{0}^{t}\!\bigl[ 
  e^{(t-\tau)A(\tau)}-e^{(t-\tau)A(t)}\big]f(t)\, d\tau\\+e^{tA(t)}f(t)-f(t)=I+II+III
\end{multline*}
where lemma \ref{lunardi} was used. The estimate for III follows. As for II one has
$$
  \| II\| \leq c\int _0^t\frac{t-\tau}{t}\frac{1}{t-\tau}\, d\tau\leq c
$$
by lemma \ref{intest}. Finally III gives
$$
  \| III\|\leq c\int _0^t\frac{1}{(t-\tau)^{1-\alpha}}\frac{1}{\tau ^\ga}\, d\tau=c\int _0^1
  \frac{1}{(1-\sigma)^{1-\alpha}}\frac{1}{\sigma ^\alpha}\, .
$$
\underline{Step 2}: Away from the origin it is again possible to argue as in the regular case.
Using 
\begin{multline*}
  u(\delta)\in D(A(\delta))\, ,\: A(\delta)u(\delta)+f(\delta)\in E^\alpha _\infty\text{ and }
  \\\| u(\delta)\| _{D(A(\delta))}+\| A(\delta)u(\delta)+f(\delta)\| _{E^\alpha _\infty}\leq c
  \| f\| _{\alpha,\rho-1}
\end{multline*}
It follows again from \cite[Theorem 4.3.1(iii)]{Lun95} that
$$
  \dot u,Au\in\opc ^\alpha\bigl([\delta,T],E_0\bigr)
$$
and that
\begin{multline*}
   \| \dot u\| _{\alpha,[\delta,T]}+\| Au\| _{\alpha,[\delta,T]}\leq c\bigl( \|f\| _{\alpha,
  [\delta,T]} +\| u(\delta)\| _{D(A(\delta))}+\| A(\delta)u(\delta)+f(\delta)\| _{E^\alpha _\infty}
  \bigr)\\\leq c\frac{1}{\delta ^\alpha}\| f\| _{\alpha,s} 
\end{multline*}
because of the embedding
$$
  \opc ^\alpha _{\alpha,\rho -1}\bigl((0,T],E_0\bigr)\hookrightarrow\opc ^\alpha _\alpha
  \bigl((0,T],E_0\bigr)
$$
and since
$$
  \| f\| _{\alpha,[\delta,T]}\leq\, c\frac{1}{\delta ^\alpha}\|f\|_{\alpha,s}
$$
as can be easily checked. The desired estimate is therefore obtained.
\end{proof}
\begin{rem}\label{limits}
It should be pointed out that, whereas conditions \eqref{hypo1}-\eqref{hypo3} are quite general, 
they, however, exclude singular families like
$$
  A(t)=\frac{A}{t^\beta}
$$
with $\beta\leq 1$ and the generator $A$ of an analytic exponentially decaying semigroup $T_A$. 
This is not 
due to \eqref{hypo1}-\eqref{hypo3} being too restrictive but rather to the fact that the regularity 
results obtained here are not valid in that case. This follows from the fact that the inequality
$$
  \| A(\tau)U(t,\tau)\| _{\mathcal{L}(E_0)}\leq\| \frac{A}{\tau ^\beta}
  T_A\bigl(\int_{\tau}^{t}\!\sigma ^{-\beta}\,d\sigma\bigr)\| _{\mathcal{L}(E_0)}\leq c
  \frac{t^\beta}{\tau ^\beta}\frac{1}{t-\tau}
$$
does not yield the needed 
$$
  \| A(\tau)U(t,\tau)\| _{\mathcal{L}(E_0)}\leq c\frac{1}{t-\tau}
$$ 
with a constant independent of $\tau$. 
If follows that this case has to be treated differently. This example also shows that condition 
\eqref{hypo2} cannot be weakened to an analogous singular H\"older condition.
\end{rem}
\section{An Example}\label{example}
In this last section an example is considered of a initial boundary value problem on a moving domain 
which undergoes an initial dimensional change. It is the latter that will eventually lead to a 
singular evolution equation of type \eqref{sacp}. Let a function $0<\inf(\varphi)\leq
\varphi\in\opbuc ^{2+\gb}(\bbr ^{n-1},\bbr)$ be given. Consider the diffusion equation
\begin{equation}\label{exeq1}
 \dot u -\triangle u =0\, ,\: (x,y)\in\bbr ^{n-1}\times\bbr\text{ with }0<y< t\varphi(x)\, ,\: t>0\, .
\end{equation}
complemented by the boundary conditions
\begin{align}\label{exeq2}
 u(t,x,0)&=g(x)\, ,\: (x,t)\in\bbr ^{n-1}\times\bbr\, ,\\\label{exeq3}
 \pd _\nu u(t,x,t\varphi(x))&=h(x)\, ,\: (x,t)\in\bbr ^{n-1}\times\bbr\, ,
\end{align}
for some given $g\in\opbuc ^{2+\gb}(\bbr ^{n-1},\bbr)$ and $h\in\opbuc ^{1+\gb}(\bbr ^{n-1},\bbr)$.
Problem \eqref{exeq1}-\eqref{exeq3} is a parabolic initial boundary value problem in the space-time 
wedge 
$$ 
 \bigcup _{t\leq 0}\{ t\}\times[0<y<t\varphi]\, .
$$
\begin{rem}\label{fbp}
For the Free Boundary Problems mentioned earlier in the paper the upper boundary of the domain would 
be given by a function $\varphi(t,x)$ which is itself an unknown of the problem and satisfies an 
additional evolution equation with initial condition $\varphi(0,\cdot)\equiv 0$, thus introducing a 
singularity into the problem via the change of variable \eqref{cov}.
\end{rem}\\
In order to apply the abstract results derived in the previous sections, the problem needs to be 
reformulated in a new set of variables
\begin{equation}\label{cov}
   (\tau ,\xi, \eta):=(t,x,\frac{y}{t\varphi(x)})\, .
\end{equation}
If one rewrites the equations in the new variables using the names of the old variables, one 
obtains
{\small\begin{multline*}
 \dot u -\triangle _x u -\frac{1+t^2|\nabla\varphi |^2}{t^2\varphi ^2}\pd _{yy}u= \frac{y}{t}\pd _yu+
 2(\frac{\nabla\varphi}{\varphi}|\pd _y \nabla u)+\frac{\varphi\triangle\varphi-|\nabla\varphi |^2}{\varphi ^2}
 \pd _y u\, , (x,y)\in S\, ,\: t>0\, ,\\
 u(t, x, 0)=g(x)\, ,\:x\in\bbr ^{n-1}\, ,\: t>0\, ,\\
 \pd _yu(t,x,1)=\frac{t\varphi}{1+t|\nabla\varphi |^2}\Bigl[ h\sqrt{1+|\nabla\varphi |^2}+
 \bigl(\nabla\varphi |\nabla _xu(t,x,1)\bigr)\Bigr]\, ,\: x\in\bbr ^{n-1}\, ,\: t>0\, ,
\end{multline*}}
for $S=\bbr ^{n-1}\times(0,1)$. As the general case is included in the forthcoming 
analysis \cite{G053} of free boundary problems, the simplifying assumption $\varphi\equiv 1$ 
is now made which leads to the simpler system
\begin{align}\label{fexeq1}
 \dot u -\triangle _x u -\frac{1}{t^2}\pd _{yy}u-\frac{y}{t}\pd _yu&=0
 \, , (x,y)\in S\, ,\: t>0\, ,\\\label{fexeq2}
 u(t, x, 0)&=g(x)\, ,\:x\in\bbr ^{n-1}\, ,\: t>0\, ,\\\label{fexeq3}
 \pd _yu(t,x,1)&=th(x)\, ,\: x\in\bbr ^{n-1}\, ,\: t>0\, ,
\end{align}
The solution of \eqref{fexeq1}-\eqref{fexeq3} can sought in the form
$$
 u(t,x,y)=v(t,x,y)+\bigl(\calr _D(t) g\bigr)(x)+t\bigl(\calr _N(t)\bigr)h(x)\, ,\: 
 (t,x,y)\in(0,\infty)\times S\, ,
$$
where the function $v$ satisfies the equation
\begin{equation}\label{mfexeq1}
 \dot v -\triangle _x v -\frac{1}{t^2}\pd _{yy}v-\frac{y}{t}\pd _yv=\bigl[\frac{y}{t}\pd _y-\pd _t\bigr]
 \bigl( \calr _D(t) g+t\calr _N(t)h\bigr)
\end{equation}
complemented with homogeneous boundary conditions, i.e. with
\begin{equation}\label{mfexeq2}
 v(t,x,0)=0=\pd _y v(t,x,1)\, .
\end{equation}
The functions $\calr _D(t) g$ and $\calr _N(t)h$ are defined as follows
\begin{gather}\label{homo1}
 \calr _D(t) g=\calf ^{-1} \frac{\cosh\bigl( t|\xi |(1-y)\bigr)}{\cosh(t|\xi |)}
 \calf g\, ,\\ \label{homo2}\calr _N(t) h=\calf ^{-1} \frac{\sinh(t|\xi |y)}{t|\xi |\cosh(t|\xi |)}\calf h\, ,
\end{gather}
and satisfy
$$
 -\triangle _x u -\frac{1}{t^2}\pd _{yy}u=0\, ,
$$
complemented with the boundary conditions
\begin{gather*}
 u(t,\cdot ,0)=g\, ,\: u_y(t,\cdot ,1)=0\text{ for }u=R_D(t)g\, ,\\
 u(t,\cdot ,0)=0\, ,\: u_y(t,\cdot ,1)=th\text{ for }u=tR_N(t)h\, ,
\end{gather*}
respectively. 
Setting $E_0=\opbuc^\gb\bigl(\bbr ^{n-1},\opc([0,1])\bigr)$ it is possible to use the results collected 
in the previous sections to obtain well-posedness of \eqref{mfexeq1}-\eqref{mfexeq2} in the 
space $\opc ^\ga _\ga\bigl((0,T],E_0\bigr)$ for $\ga\in(0,1)$.
\begin{lem}
The function $f:=\bigl[\frac{y}{t}\pd _y-\pd _t\bigr]\bigl[ \calr _D(t) g+t\calr _N(t)h\bigr]$ satisfies
$$
 f\in\opc ^\ga_\ga\bigl((0,T],E_0\bigr)
$$
\end{lem}
\begin{proof}
A direct computation shows that
\begin{multline*}
 \calf f=\bigl[ (1-y)\sinh\bigl( t|\xi |(1-y)\bigr) +\tanh(t|\xi |) \cosh\bigl( t|\xi |(1-y)\bigr)\bigr]
 \frac{1}{\cosh(t|\xi |)}|\xi |\,\calf g +\\ - \bigl[ y\cosh(t|\xi |y) -\tanh(t|\xi |) \sinh(t|\xi |y) \bigr] 
 \frac{1}{\cosh(t|\xi |)}\calf h+\\-y\frac{\sinh\bigl( t|\xi |(1-y)\bigr)}{\cosh(t|\xi |)}|\xi |\,\calf g+
 y\frac{\cosh(t|\xi |y)}{\cosh(t|\xi |)}\calf h
\end{multline*}
The claim then follows from \cite[Lemmata 2.2, 2.5, 2.6]{Gui99} combined with the operator-valued 
Fourier multiplier theorem \cite[Theorem 6.2]{Ama97}. The fact that only the case $n=2$ is considered 
in \cite{Gui99} is an immaterial difference, since the results involving $n$ remain valid with only 
obvious modifications for $n>2$.
\end{proof}

It is well-known that 
\begin{equation}\label{op1}
 B=\triangle _x:\dom(B)\subset\opbuc^\gb(\bbr ^{n-1})\to\opbuc^\gb(\bbr ^{n-1})
\end{equation}
generates a bounded analytic semigroup on $\opbuc^\gb\bigl(\bbr,\opc([0,1])\bigr)$ and it is easy the check
by standard elliptic theory (see \cite[Chap. 3]{Lun95} for instance) that 
\begin{equation}\label{op2}
 C=\pd _{yy}+ty\pd _y:\dom(C)\subset\opc([0,1])\to\opc([0,1])
\end{equation}
(with boundary conditions) generates an exponentially decaying analytic semigroup on $\opbuc^\gb\bigl(
\bbr ^{n-1},\opc([0,1])\bigr)$ for any fixed $t\in[0,T]$. In either case the missing variables can be 
considered as parameters. The sum $A (t)=B+\frac{1}{t^2}C$ then generates an exponentially decaying 
analytic semigroup on $E_0$ with constant domain of definition, provided $t>0$. 
In order to verify that all conditions \eqref{hypo1}-\eqref{hypo3} are satisfied we still need 
to estimate
$$
 \bigl[A(t)-A(s)\bigr]\bigl(-A (\tau)\bigr)^{-q}=\Bigl[\bigl( \frac{1}{s^2}-\frac{1}{t^2}\bigr)
 \pd _{yy}+\bigl( \frac{1}{s}-\frac{1}{t}\bigr)y\pd _y\Bigr]\Bigl[-\triangle _x-\frac{1}{\tau ^2}
 \pd _{yy}-\frac{1}{\tau}y\pd _y\bigr]^{-q}
$$
for $0<\tau\leq s\leq t\leq T$ and $q=1,p$. Since the estimate is only needed in the base 
space $E_0$, it causes no problem that
$$
 \Bigl[ -\tau ^2\triangle _x-\pd _{yy}-\tau y\pd _y\bigr]^{-q}
$$
looses its $x$-regularizing effect as $\tau\to 0$ provided $\pd _{yy}$ (with boundary conditions) is 
invertible, as it is. This leads to the estimates
$$
  \| \bigl[A(t)-A(s)\bigr]\bigl(-A (\tau)\bigr)^{-q}\| _{\call(E_0)}\leq \,c 
  (t-s)/t^{\frac{p-q}{p-1}}\, ,\: q=1,p\, .
$$
Theorem \ref{max_reg_sing} can now be safely applied to obtain the following result.
\begin{thm}
Assume that $g\in\opbuc ^{2+\gb}(\bbr ^{n-1},\bbr)$ and $h\in\opbuc ^{1+\gb}(\bbr ^{n-1},\bbr)$ 
for $\gb\in(0,1)$. Then problem \eqref{fexeq1}-\eqref{fexeq3} possesses a unique solution with
$$
 u,\dot u,A u\in\opc ^\ga _\ga\Bigl((0,T],\opbuc ^\gb\bigl(\bbr ^{n-1},\opc([0,1])\bigr)\Bigr)\, .
$$ 
It is given by
$$
 u(t)=\int _0^tU _A(t,\tau)f(\tau)\, d\tau + \calr _D(t)g +t\calr _N(t)h\, ,
$$
where $U$ is the evolution operator generated by the singular family $A (t)=B+\frac{1}{t^2}C$ 
with $B$ and $C$ as in \eqref{op1}-\eqref{op2} and 
the operators $\calr _D$ and $\calr _N$ are as defined in \eqref{homo1}-\eqref{homo2}.
\end{thm}
\bibliography{../../lite}
\end{document}